\documentclass[11pt]{article}
\usepackage{amssymb}
\usepackage{amsfonts}
\usepackage{mathrsfs}
\makeatletter
       \@addtoreset{equation}{section}\makeatother
\newtheorem{Definition}{Definition}[section]
\newtheorem{Example}{Example}[section]
\newtheorem{Theorem}{Theorem}[section]
\newtheorem{Lemma}[Theorem]{Lemma}
\newtheorem{Corollary}[Theorem]{Corollary}

\def\calL{\mathcal{L}}
\def\calH{\mathcal{H}}
\def\bbR{\mathbb{R}}
\def\bbZ{\mathbb{Z}}
\def\bbJ{\mathbb{J}}
\def\bbK{\mathbb{K}}
\def\LTR{L^2(\bbR)}
\begin{document}
\parindent 16pt
\parskip 5pt
\baselineskip 15pt
\title{G-frames and G-Riesz Bases\thanks{This work was supported partially  by
the National Natural Science Foundation of China (10201014 and
60472042),
the Program for New Century Excellent Talents in Universities,
and the Research Fund for the Doctoral Program of
Higher Education.
}}

\author{Wenchang  Sun\\
\mbox{}\\
  Department of Mathematics and LPMC,  Nankai University,\\
      Tianjin 300071, China\\
      Email:\, sunwch@nankai.edu.cn}

\date{June 28, 2005}
\maketitle
\textbf{Abstract}
G-frames are generalized frames  which include
ordinary frames,
bounded invertible linear operators,
as well as
many recent generalizations of frames,
  e.g.,
        bounded quasi-projectors
and
        frames of subspaces.
G-frames are natural generalizations of frames
and
provide  more choices on  analyzing   functions
from frame expansion coefficients.
We give  characterizations of g-frames
and prove that g-frames share many useful properties with   frames.
We also give generalized version of Riesz bases and orthonormal bases.
As an application, we get atomic resolutions for bounded linear operators.



\section{Introduction}

Frames were first introduced in 1952 by Duffin and Schaeffer
\cite{DS},
reintroduced in 1986 by Daubechies, Grossman, and Meyer \cite{CGY},
and
popularized from then on.
Frames have many nice properties which
make them very useful in the characterization of function spaces,
signal processing and  many other fields. We refer to
 \cite{C6,D2,FeiS,Gr3,HL,HW,Y} for an
introduction to the frame theory and its applications.
One of the main virtues of frames is that,
given a frame, we can get properties of
a function and reconstruct it only from the frame
coefficients, a sequence of complex numbers.
For example, let $\{a^{j/2}\psi_{\ell}(a^j\cdot-bk):\,1\le\ell\le r, j,k\in\bbZ\}$
be a multi-wavelet frame for
$\LTR$. Then every $f\in \LTR$ can be reconstructed by the sequence
$\{\langle f, a^{j/2}\psi_{\ell}(a^j\cdot-bk)\rangle:\,1\le\ell\le r, j,k\in\bbZ\}$
which satisfying
\[
A\|f\|_2^2 \le
\sum_{1\le\ell\le r}\sum_{j,k\in\bbZ} |\langle f, a^{j/2}\psi_{\ell}(a^j\cdot-bk)\rangle|^2
\le
B\|f\|_2^2
\]
for some positive constants $A$ and $B$.
Put
\[c_{j,k}(f) = (\langle f, a^{j/2}\psi_{1}(a^j\cdot-bk)\rangle, \cdots,
\langle f, a^{j/2}\psi_{r}(a^j\cdot-bk)\rangle)^T\in \mathbb C^r.
\]
Then the  above
inequalities  turn out to be
\[
A\|f\|_2^2 \le
\sum_{j,k\in\bbZ} \|c_{j,k}(f)\|^2
\le
B\|f\|_2^2.
\]
This prompts us to give the following generalization of frames.

Throughout this paper,
$\mathcal{U}$ and $\mathcal{V}$ are two Hilbert spaces
and $\{\mathcal{V}_j:\,j\in\bbJ\}\subset \mathcal{V}$
is a sequence of Hilbert spaces,
where $\bbJ$ is a subset of $\bbZ$.
$\calL(\mathcal{U},\mathcal{V}_j)$ is the collection of all bounded
linear operators from $\mathcal{U}$ into $\mathcal{V}_j$.

Note that for any sequence
$\{\mathcal{V}_j:\,j\in\bbJ\}$
 of Hilbert spaces, we can always find a big space $\mathcal{V}$ to contain all the $\mathcal{V}_j$
 by setting
 $\mathcal{V} = \bigoplus_{j\in\bbJ} \mathcal{V}_j$.
\begin{Definition}
\label{def:g-frame}
We call a sequence $\{\Lambda_j\in\calL(\mathcal{U},\mathcal{V}_j) :\,j\in\bbJ\}$
  a generalized frame, or simply a g-frame,
for $\mathcal{U}$ with respect to $\{\mathcal{V}_j:\,j\in\bbJ\}$ if there are two positive
constants $A$ and $B$ such that
\begin{equation}\label{eq:frame}
  A\|f\|^2 \le \sum_{j\in\bbJ} \|\Lambda_j f\|^2 \le B\|f\|^2,
  \qquad \forall f\in \mathcal{U}.
\end{equation}
We call $A$ and $B$ the lower and upper frame bounds, respectively.

We call $\{\Lambda_j:\,j\in\bbJ\}$ a tight g-frame if $A=B$.

We call $\{\Lambda_j:\,j\in\bbJ\}$ an exact g-frame
 if it ceases to be a g-frame whenever anyone of its elements is removed.

We say simply a g-frame for $\mathcal{U}$ whenever the space sequence
$\{\mathcal{V}_j:\,j\in\bbJ\}$ is clear.

We say also a g-frame for $\mathcal{U}$ with respect to $\mathcal{V}$
whenever $\mathcal{V}_j=\mathcal{V}$, $\forall j\in\bbJ$.
\end{Definition}


We observe that various generalizations of frames have been proposed recently.
For example,
bounded quasi-projectors \cite{For1,For2},
frames of subspaces\cite{AK,CK},
pseudo-frames\cite{LO},
oblique frames\cite{CE,Eld},
and outer frames\cite{ACM}.
All of these generalizations are proved to be useful in many applications.
Here we point out that they can be regarded as special cases of g-frames (see examples below)
and many basic properties can be derived within this more general context.

While we were preparing  this paper we learned
that another generalization of frames in the context of numerical analysis,
called stable space splittings,
have been studied in \cite{Osw1,Osw2}.
We prove at the end of Section 3 that they are equivalent to g-frames.
We point out that the approaches are quite different from each others.
In particular, the adjoint operators of $\Lambda_j$ are used in the definition
of
stable space splittings.
Moreover,  we give a characterization of g-frames
and
 studied g-Riesz bases and g-orthonormal bases.

\begin{Example} \label{ex:ex1}
Let $\calH$ be a separable Hilbert space and $\{f_j:\,j\in\bbJ\}$ be a
 frame for $\calH$.
Let $\Lambda_{f_j}$ be the functional induced by $f_j$, i.e.,
\[
      \Lambda_{f_j} f = \langle f, f_j\rangle, \qquad \forall f\in \calH.
\]
It is easy to check that
 $\{\Lambda_{f_j}:\,j\in\bbJ\}$ is a g-frame for $\calH$ with respect to $\mathbb{C}$.
\end{Example}

By Riesz Representation Theorem,
to every functional $\Lambda \in\mathcal{L}(\mathcal{U},\mathbb{C})$,
one can find some $\varphi\in \mathcal{U}$ such that
$\Lambda f = \langle f, \varphi\rangle$, $\forall f\in\mathcal{U}$.
Hence we have the following.
\begin{Lemma}
\label{Lm:o-frame}
A frame is equivalent to
a g-frame whenever $\mathcal{V}_j=\mathbb{C}$, $j\in\bbJ$.
\end{Lemma}

\begin{Example}
\label{ex:ex2}
Pseudo-frames (Li and Ogawa\cite{LO}), or similar,
oblique frames(Christensen and Eldar \cite{CE,Eld})
or
outer frames
 (Aldroubi, Cabrelli, and Molter \cite{ACM})
 are studied recently in literature.
Here we point out that they are   a class of g-frames.
\end{Example}

Let $\calH_0$ be a closed subspace of $\calH$.
Let $\{f_j:\,j\in\bbJ\}\subset \calH$ be a Bessel sequence in $\calH_0$
and $\{\tilde f_j:\,j\in\bbJ\}\subset \calH$ be a Bessel sequence in $\calH$.
Recall that  $\{f_j:\,j\in\bbJ\}$
is said to be a   pseudo-frame
for $\calH_0$ with respect to
$\{\tilde f_j:\,j\in\bbJ\}$ \cite[Definition 1]{LO}
if
\[
   f = \sum_{j\in\bbJ} \langle f, f_j\rangle \tilde f_j, \qquad \forall f\in \calH_0.
\]
Since both
$\{f_j:\,j\in\bbJ\}$ and
$\{\tilde f_j:\,j\in\bbJ\}$ are Bessel sequences in $\calH_0$, it is easy to check
from the above equation  that we can find some
constants $A,B>0$ such that
$A\|f\|^2 \le   \sum_{j\in\bbJ}
   | \langle  f, f_j \rangle |^2 \le B\|f\|^2$, $\forall f\in  \calH_0$.
Let $\Lambda_{f_j}$  be the  functional induced by $f_j$, $j\in\bbJ$.
Then we have
\[
A\|f\|^2 \le   \sum_{j\in\bbJ}
   | \Lambda_{f_j} f|^2 \le B\|f\|^2, \qquad \forall f\in
       \calH_0.
\]
In other words,
$\{\Lambda_{f_j}:\,j\in\bbJ\}$ is a g-frame for $\calH_0$ with respect to
$\mathbb{C}$.

\begin{Example}
\label{ex:quasi-projectors}
Bounded quasi-projectors (Fornasier \cite{For1,For2}).
\end{Example}

It was shown in \cite[Lemma 1]{For1} that if a system of
bounded quasi-projectors $\{P_j:\,j\in\bbJ\}$ is self-adjoint
and compatible with the canonical projections
(see \cite{For1,For2} for details),
then for any $f\in\mathcal{H}$,
\[
  A \|f\|^2 \le \sum_{j\in\bbJ} \|P_j f\|^2 \le B\|f\|^2.
\]
In this case,
 $\{P_j:\,j\in\bbJ\}$ is a g-frame for $\mathcal{H}$ with respect to $\mathcal{H}$.

\begin{Example}
\label{ex:CK}
Frames of subspaces
(Casazza and Kutyniok \cite{CK}
and  Asgari and Khosravi\cite{AK}).
\end{Example}

Let $\{W_j:\,j\in\bbJ\}$
be a sequence of subspaces of $\calH$ and
$P_{W_j}$ be the orthogonal projection on $W_j$.
$\{W_j:\,j\in\bbJ\}$ is called a frame of subspaces if
there exist positive constants $A$ and $B$ such that
\[
A\|f\|^2\le   \sum_{j\in\bbJ} \|P_{W_j} f\|^2 \le B\|f\|^2,
\qquad \forall f\in\calH.
\]
Obviously, a frame of subspaces is  a g-frame
for $\mathcal{H}$ with respect to $\{W_j:\,j\in\bbJ\}$.

\begin{Example}
\label{Ex:DFG}
Time-frequency localization operators
(D\"orfler, Feichtinger and Gr\"ochenig \cite{DFG}).
\end{Example}

For $f,g\in L^{2}(\bbR^d)$, define the windowed Fourier transform of $f$ with respect to $g$ by
\[
   (V_g f)(t,\omega) = \int_{\bbR^d} f(x) \overline{g(x-t)} e^{-i2\pi x\omega} d x.
\]
Let $S_0(\bbR^d):= \{g\in L^2(\bbR^d):\, V_g g \in L^1(\bbR^{2d})\}$ be the Feichtinger algebra.
Take some $\varphi\in S_0(\bbR^d)$ with  $\|\varphi\|_2=1$.
Let $\sigma$ be a bounded function on $\bbR^{2d}$
with compact support and $\sigma(x)\ge 0$.
Define the time-frequency localization operator $H_{\sigma}$
corresponding to $\sigma$ and $\varphi$ by
$
   H_{\sigma}f = V_\varphi^* \sigma V_\varphi f.
$
If $\sigma\in S_0(\bbR^{2d})$ and
\[
    C_1\le \sum_{k\in\bbZ^{2d}} \sigma(x-k) \le C_2,
\]
for some constants $C_1, C_2>0$, then it is shown in \cite{DFG}
that one can find some constants $A, B>0$ such that
\[
  A\|f\|_2^2 \le \sum_{k\in\bbZ^{2d}} \| H_{\sigma(\cdot-k)} f\|_2^2
  \le B\|f\|_2^2, \qquad \forall f\in L^2(\bbR^d).
\]
Hence
$\{H_{\sigma(\cdot-k)}:\, k\in\bbZ^{2d}\}$ is a g-frame for $L^2(\bbR^d)$
with respect to
$L^2(\bbR^d)$. We refer to \cite{DFG} for details.

\begin{Example}
\label{ex:op}
Every bounded invertible linear operator itself forms
a g-frame.
\end{Example}

We see from the above examples that
g-frames are natural generalizations of frames
and  provide more choices on  analyzing   functions
from frame expansion coefficients
.
In the following sections we first
study g-frame operators and get the dual g-frames,
then give  definitions of
g-Riesz bases and g-orthonormal bases
and present characterizations of generalized frames and bases.
As an application of g-frames,
we get atomic resolutions of bounded linear operators.

\section{G-frame operators and dual g-frames}
\def\seqV{\{\mathcal{V}_j:\,j\in\bbJ\}}

Let $\{\Lambda_j:\,j\in\bbJ\}$ be a   g-frame for
$\mathcal{U}$ with respect to $\seqV$. Define the  g-frame operator $S$ as follows:
\begin{equation}\label{eq:fs}
   S f = \sum_{j\in\bbJ} \Lambda_j^*  \Lambda_j f, \quad \forall f \in \mathcal{U},
\end{equation}
where $\Lambda_j^*$ is the adjoint operator of $\Lambda_j$.
First of all, $S$ is well defined on $\mathcal{U}$.
To see this, let $n_1<n_2$ be integers. Then
we have
\begin{eqnarray*}
   \left\| \sum_{j=n_1}^{n_2} \Lambda_j^* \Lambda_j f\right\|
   &=& \sup_{h\in\mathcal{U}, \|h\|=1}
      \left| \left\langle\sum_{j=n_1}^{n_2}
        \Lambda_j^* \Lambda_j f , h\right\rangle\right| \\
   &=&\sup_{\|h\|=1}
         \left| \sum_{j=n_1}^{n_2}\langle  \Lambda_j f , \Lambda_j h\rangle\right| \\
   &\le&\sup_{\|h\|=1}   \left(\sum_{j=n_1}^{n_2}\| \Lambda_j f\|^2 \right)^{1/2}
           \cdot
            \left(\sum_{j=n_1}^{n_2}\| \Lambda_j h\|^2 \right)^{1/2}\\
   &\le& B^{1/2}\left(\sum_{j=n_1}^{n_2}\| \Lambda_j f\|^2 \right)^{1/2}.
\end{eqnarray*}
Now we see from (\ref{eq:frame}) that
the series in (\ref{eq:fs}) are convergent.
Therefore, $Sf$ is well defined for any $f\in \mathcal{U}$.

On the other hand, it is easy to check that for any $f_1, f_2\in \mathcal{U}$,
\[
\langle Sf_1, f_2\rangle
 = \sum_{j\in\bbJ} \langle \Lambda_j^*  \Lambda_j f_1, f_2\rangle
 = \sum_{j\in\bbJ} \langle  f_1, \Lambda_j^*  \Lambda_j f_2\rangle
=\langle f_1, S f_2\rangle
 \]
and therefore,
\[
  \|S\| = \sup_{\|f\|=1} \langle Sf, f\rangle
    = \sup_{\|f\|=1} \sum_{j\in\bbJ} \|\Lambda_j f\|^2 \le B.
\]
Hence $S$ is a bounded self-adjoint operator.

Since $A\|f\|^2  \le \langle Sf, f\rangle \le \|Sf\|\cdot \|f\|$, we have
\[
         \|Sf\| \ge A \|f\|,
\]
which implies that $S$ is injective and
$S\mathcal{U}$ is closed in $\mathcal{U}$.
Let $f_2\in \mathcal{U}$ be such that $\langle Sf_1, f_2\rangle =0$ for any $f_1\in \mathcal{U}$.
Then we have $\langle f_1, S f_2\rangle =0$, $\forall f_1\in \mathcal{U}$. This implies that
$S f_2=0$ and therefore $f_2=0$. Hence $S\mathcal{U}=\mathcal{U}$.
Consequently, $S$ is invertible and
\[
        \|S^{-1} \| \le \frac{1}{A}.
\]
For any $f\in \mathcal{U}$, we have
\[
     f= SS^{-1}f = S^{-1}Sf= \sum_{j\in\bbJ} \Lambda_j^* \Lambda_j S^{-1} f
 = \sum_{j\in\bbJ} S^{-1}\Lambda_j^* \Lambda_j  f.
\]
Let
\(
       \tilde{\Lambda}_j = \Lambda_j S^{-1}.
\)
Then the above equalities become
\begin{equation} \label{Eq:expan}
     f=  \sum_{j\in\bbJ} \Lambda_j^* \tilde\Lambda_j  f
     =
      \sum_{j\in\bbJ} \tilde\Lambda_j^* \Lambda_j  f.
\end{equation}
We now prove that $\{\tilde\Lambda_j:\,j\in\bbJ\}$ is also a g-frame for
$\mathcal{U}$
with respect to $\seqV$.

In fact, for any $f\in \mathcal{U}$, we have
\begin{eqnarray*}
\sum_{j\in\bbJ} \|\tilde\Lambda_j f\|^2
&=&\sum_{j\in\bbJ} \| \Lambda_j S^{-1} f\|^2 \\
&=&\sum_{j\in\bbJ} \langle \Lambda_j S^{-1} f, \Lambda_j S^{-1} f\rangle \\
&=&\sum_{j\in\bbJ} \langle \Lambda_j^*\Lambda_j S^{-1} f,  S^{-1} f\rangle \\
&=& \langle S S^{-1}f, S^{-1}f\rangle \\
&=& \langle  f, S^{-1}f\rangle\\
&\le& \frac{1}{A}\|f\|^2.
\end{eqnarray*}
On the other hand, since
\begin{eqnarray*}
\|f\|^2
&=& \sum_{j\in\bbJ} \langle \tilde\Lambda_j^*\Lambda_j f, f\rangle\\
&=& \sum_{j\in\bbJ} \langle \Lambda_j f, \tilde\Lambda_j  f\rangle\\
&\le& \left(\sum_{j\in\bbJ} \|  \Lambda_j f \|^2 \right)^{1/2}
      \cdot \left(\sum_{j\in\bbJ} \|  \tilde\Lambda_j f \|^2 \right)^{1/2}
        \\
&\le& B^{1/2} \|f\| \left(\sum_{j\in\bbJ} \|  \tilde\Lambda_j f \|^2 \right)^{1/2},
\end{eqnarray*}
we have
\[
      \sum_{j\in\bbJ} \|  \tilde\Lambda_j f \|^2
      \ge  \frac{1}{B}\|f\|^2.
\]
Hence,
$\{\tilde\Lambda_j:\,j\in\bbJ\}$ is a g-frame for
$\mathcal{U}$ with frame bounds $1/B$ and $1/A$.
We call it the (canonical) dual g-frame of
$\{\Lambda_j:\,j\in\bbJ\}$.

Let
$\tilde{S}$
be  the g-frame operator associated with
$\{\tilde\Lambda_j:\,j\in\bbJ\}$.
Then we have
\begin{eqnarray*}
   S\tilde{S} f
   &=& \sum_{j\in\bbJ} S \tilde\Lambda_j^* \tilde\Lambda_j f
   = \sum_{j\in\bbJ} S \,\, S^{-1} \Lambda_j^*\,\, \Lambda_j S^{-1}\ f  \\
   &=& \sum_{j\in\bbJ}   \Lambda_j^*\Lambda_j S^{-1} f
   = SS^{-1}f
   = f, \qquad \forall f\in \mathcal{U}.
\end{eqnarray*}
Hence $\tilde{S} = S^{-1}$ and $\tilde\Lambda_j \tilde{S}^{-1}= \Lambda_j S^{-1} S = \Lambda_j$.
In other words,
$\{\Lambda_j:\,j\in\bbJ\}$ and
$\{\tilde\Lambda_j:\,j\in\bbJ\}$ are dual g-frames with respect to each other.

\textbf{Remark}.\,\,
We see from the above arguments that g-frames behave very similarly to frames.
For example,
we can always get a tight g-frame from any g-frame
$\{\Lambda_j:\,j\in\bbJ\}$. In fact, put
\[
     Q_j = \Lambda_jS^{-1/2}.
\]
It is easy to check that $\{Q_j:\,j\in\bbJ\}$ is a tight g-frame with the frame bound $1$.

Moreover,
the canonical dual g-frames give rise to expansion coefficients with the minimal
norm.
\begin{Lemma} \label{Lm:g-frame1}
Let $\{\Lambda_j:\,j\in\bbJ\}$ be a g-frame for $\mathcal{U}$ with respect to
$\seqV$ and $\tilde{\Lambda}_j = \Lambda_j S^{-1}$.
Then for any $g_j\in \mathcal{V}_j$ satisfying
$f = \sum_{j\in\bbJ} \Lambda_j^* g_j$, we have
\[
     \sum_{j\in\bbJ} \|g_j\|^2 = \sum_{j\in\bbJ} \|\tilde\Lambda_j f\|^2
           +\sum_{j\in\bbJ} \|g_j -\tilde\Lambda_j f\|^2.
\]
\end{Lemma}

\begin{proof}
It is easy to check that
\begin{eqnarray*}
  \sum_{j\in\bbJ} \|\tilde\Lambda_j f\|^2
   &=& \sum_{j\in\bbJ} \langle \tilde\Lambda_j f, \Lambda_j S^{-1} f\rangle \\
   &=& \sum_{j\in\bbJ} \langle \Lambda_j^*\tilde\Lambda_j f,  S^{-1} f\rangle \\
   &=& \sum_{j\in\bbJ} \langle \Lambda_j^* g_j,  S^{-1} f\rangle \\
   &=& \sum_{j\in\bbJ} \langle  g_j,  \Lambda_j S^{-1} f\rangle \\
   &=& \sum_{j\in\bbJ} \langle  g_j,  \tilde\Lambda_j f\rangle, \qquad \forall f\in \mathcal{U}.
\end{eqnarray*}
Now the conclusion follows.
\end{proof}

In example \ref{ex:ex1}, we show that every   frame
$\{f_j:\, j\in\bbJ\}$ for $\calH$
induces a g-frame
$\{\Lambda_{f_j}:\, j\in\bbJ\}$ for $\calH$ with respect to $\mathbb{C}$ via the
induced functionals $\Lambda_{f_j}$.

Let $\{\tilde{f_j}:\,j\in\bbJ\}$ be the  canonical dual   frame of
$\{f_j:\, j\in\bbJ\}$.
We conclude that
$\{\Lambda_{\tilde{f_j}}:\,j\in\bbJ\}$
is the canonical dual g-frame of
$\{\Lambda_{f_j}:\, j\in\bbJ\}$.

In fact, it is easy to see that
$\Lambda_{f_j}^* c = c f_j$ for any $c\in \mathbb{C}$,
which implies that
the corresponding
g-frame operator and frame operator  are the same.
Consequently,
\[
   \Lambda_{f_j}S^{-1} f = \langle S^{-1}f, f_j\rangle
    = \langle f, S^{-1} f_j\rangle
    =\langle f, \tilde{f}_j\rangle
    =\Lambda_{\tilde{f_j}} f, \qquad \forall f\in\calH.
\]
Hence
\(
\Lambda_{\tilde{f_j}} = \Lambda_{f_j}S^{-1}.
\)
In other words,
$\{\Lambda_{\tilde{f_j}}:\,j\in\bbJ\}$
is the  dual g-frame of
$\{\Lambda_{f_j}:\, j\in\bbJ\}$.

\section{Generalized Bessel sequences, Riesz bases and orthonormal bases}

Similarly to generalized frames,
we can define generalized Bessel sequences, Riesz bases, and orthonormal bases.

\begin{Definition}
\label{def:g-riesz}
Let $\Lambda_j\in\calL(\mathcal{U},\mathcal{V}_j)$, $j\in\bbJ$.

\begin{enumerate}
\item
If  the right-hand inequality of (\ref{eq:frame})
holds, then we say that
$\{\Lambda_j:\,j\in\bbJ\}$ is  a g-Bessel sequence for
$\mathcal{U}$ with respect to $\seqV$.

\item
    If $\{f:\,  \Lambda_j f =0, j\in \bbJ \}=\{0\}$,
    then we say that
     $\{\Lambda_j:\,j\in\bbJ\}$ is g-complete.

\item
   If
    $\{\Lambda_j:\,j\in\bbJ\}$ is g-complete and
     there are positive constants
  $A$ and $B$ such that
  for any finite subset $\bbJ_1\subset \bbJ$ and $g_j\in \mathcal{V}_j, j\in\bbJ_1$,
\begin{equation}\label{eq:g-riesz}
  A\sum_{j\in \bbJ_1}\|g_j\|^2
   \le \left\| \sum_{j\in\bbJ_1}  \Lambda_j^* g_j\right\|^2
 \le   B\sum_{j\in \bbJ_1}\|g_j\|^2,
\end{equation}
then we say that
    $\{\Lambda_j:\,j\in\bbJ\}$ is  a g-Riesz basis for $\mathcal{U}$ with respect to
  $\seqV$.

\item
we say
 $\{\Lambda_j:\,j\in\bbJ\}$ is  a g-orthonormal basis for $\mathcal{U}$ with respect to
  $\seqV$ if it satisfy the following.
\begin{eqnarray}
&& \langle \Lambda_{j_1}^* g_{j_1}, \Lambda_{j_2}^* g_{j_2}\rangle =
     \delta_{j_1, j_2} \langle g_{j_1}, g_{j_2}\rangle,
      \quad
      \forall  j_1, j_2\in\bbJ,  g_{j_1}\in\mathcal{V}_{j_1},
           g_{j_2}\in\mathcal{V}_{j_2}, \label{eq:orth:a} \\
&& \sum_{j\in\bbJ}\|\Lambda_j f\|^2=\|f\|^2, \qquad \forall f\in\mathcal{U}.\label{eq:orth:c}
\end{eqnarray}
\end{enumerate}
\end{Definition}

\begin{Example}
As in Example \ref{ex:ex1}, the induced functionals of any
  Bessel sequence(resp.   Riesz basis,   orthonormal basis)
form
a g-Bessel sequence(resp. g-Riesz basis, g-orthonormal basis).
\end{Example}

\begin{Example}
The sequence containing only the identity mapping $\{I_{\mathcal{U}}\}$  is a
 g-Bessel sequence, g-Riesz basis, and a g-orthonormal basis
 for
 $\mathcal{U}$ with respect to
  $\mathcal{U}$.
\end{Example}

\begin{Example}
Let $(X, \mathscr{B}, m)$ be a measure space
and $\{X_j:\,j\in\bbJ\}$ be a sequence of measurable sets.
Let $\Lambda_j$ be the orthonormal projection from
$L^2(X)$ onto $L^2(X_j)$, i.e.,
$
     \Lambda_j f = f \cdot \chi^{}_{X_j}.
$
Then we have
\begin{enumerate}
\item
$\{\Lambda_j:\,j\in\bbJ\}$ is a g-frame for $L^2(X)$ with respect to $\{L^2(X_j):\,j\in\bbJ\}$
if and only if
$\bigcup_{j\in\bbJ} X_j = X$
and $\sup_{j\in \bbJ} \#\{j':\, m(X_j\cap X_{j'})>0\}  <+\infty$.

\item $\{\Lambda_j:\,j\in\bbJ\}$ is a g-Riesz basis for $L^2(X)$
with respect to $\{L^2(X_j):\,j\in\bbJ\}$ if and only if
$\bigcup_{j\in\bbJ} X_j = X$
and $m(X_j\cap X_{j'})=0$, $ j \ne  j'$.
If it is the case, it is also a g-orthonormal basis.
\end{enumerate}
\end{Example}

\subsection{Characterizations of g-frames, g-Riesz bases and g-orthonormal bases}

Let $\Lambda_j\in
\mathcal{L}(\mathcal{U}, \mathcal{V}_j)$.
We do not have other assumptions on $\Lambda_j$ at the moment.
Suppose that
$\{e_{j,k}:\,k\in\bbK_j\}$
is an orthonormal basis for  $\mathcal{V}_j$,
where $\bbK_j$ is a  subset of $\bbZ$, $j\in\bbJ$.
Then
\[
    f \mapsto \langle \Lambda_j f, e_{j,k}\rangle
\]
defines a bounded linear functional on $\mathcal{U}$.
Consequently,    we can find some $u_{j,k}\in\mathcal{U}$
such that
\begin{equation}\label{eq:ujk}
    \langle f, u_{j,k}\rangle =\langle \Lambda_j f, e_{j,k}\rangle,
    \qquad \forall f\in \mathcal{U}.
\end{equation}
Hence
\begin{equation}\label{eq:rep:phij}
   \Lambda_j f = \sum_{k\in\bbK_j} \langle f, u_{j,k}\rangle e_{j,k},
   \qquad \forall f\in\mathcal{U}.
\end{equation}
Since $\sum_{k\in\bbK_j} |\langle f, u_{j,k}\rangle|^2
=\|\Lambda_j f\|^2 \le \|\Lambda_j\|^2\cdot \|f\|^2$,
$\{u_{j,k}:\,k\in\bbK_j\}$ is a    Bessel sequence for $\mathcal{U}$.
It follows that
for any $f\in\mathcal{U}$ and $g\in\mathcal{V}_j$,
\[
 \langle f, \Lambda_j^* g\rangle
= \langle \Lambda_jf, g\rangle
= \sum_{k\in\bbK_j} \langle f, u_{j,k}\rangle\cdot  \langle e_{j,k}, g\rangle
=\left\langle f,  \sum_{k\in\bbK_j}  \langle g,  e_{j,k} \rangle u_{j,k}\right\rangle.
\]
Hence
\begin{equation}\label{eq:rep:phi*}
  \Lambda_j^* g = \sum_{k\in\bbK_j} \langle g,  e_{j,k} \rangle u_{j,k},
  \qquad \forall g\in\mathcal{V}_j.
\end{equation}
In particular,
\begin{equation}\label{eq:rep:ujk}
  u_{j,k} = \Lambda_j^*    e_{j,k}, \qquad j\in\bbJ, k\in\bbK_j.
\end{equation}
We call  $\{u_{j,k}:\,j\in\bbJ, k\in\bbK_j\}$
the sequence  induced by $\{\Lambda_j:\,j\in\bbJ\}$
with respect to $\{e_{j,k}:\,j\in\bbJ, k\in\bbK_j\}$.

With above representations of $\Lambda_j$ and $\Lambda_j^*$,
we get characterizations of
generalized frames, Riesz bases, and orthonormal bases.
\begin{Theorem}
\label{thm:equivalence}
Let $\Lambda_j\in \mathcal{L}(\mathcal{U}, \mathcal{V}_j)$
and $u_{j,k}$ be defined as in (\ref{eq:rep:ujk}).
Then we have the followings.
\begin{enumerate}
\item
$\{\Lambda_j:\,j\in\bbJ\}$ is a g-frame
(resp. g-Bessel sequence,  tight g-frame, g-Riesz basis, g-orthonormal basis)
 for $\mathcal{U}$ if and only if
$\{u_{j,k}:\,j\in\bbJ,k\in\bbK_j\}$ is a    frame
(resp.   Bessel sequence,   tight   frame,
  Riesz basis,  orthonormal basis
) for $\mathcal{U}$.

\item
If
$\{\Lambda_j:\,j\in\bbJ\}$ is a g-frame,
then
\[
       \sum_{j\in\bbJ} \dim \mathcal{V}_j \ge \dim \mathcal{U}
\]
and the equality holds whenever
$\{\Lambda_j:\,j\in\bbJ\}$ is a g-Riesz basis.

\item
Moreover,
the g-frame operator for $\{\Lambda_j:\,j\in\bbJ\}$
coincides with the   frame operator for
$\{u_{j,k}:\,j\in\bbJ,k\in\bbK_j\}$.

\item
 Furthermore,
$\{\Lambda_j:\,j\in\bbJ\}$
and $\{\tilde\Lambda_j:\,j\in\bbJ\}$
are a pair of (canonical) dual g-frames
if and only if the induced sequences are a pair of
(canonical) dual   frames.
\end{enumerate}
\end{Theorem}

\begin{proof}
(i).\,\, We see from (\ref{eq:rep:phij}) that
\[
   \sum_{j\in\bbJ} \|  \Lambda_j f\|^2
  =
         \sum_{j\in\bbJ} \sum_{k\in\bbK_j} |\langle f, u_{j,k}\rangle|^2,
         \qquad \forall f\in\mathcal{U}.
\]
Hence
$\{\Lambda_j:\,j\in\bbJ\}$ is a g-frame (resp. g-Bessel sequence,
  tight g-frame)
for $\mathcal{U}$ if and only if
$\{u_{j,k}:\,j\in\bbJ,k\in\bbK_j\}$ is a    frame
(resp.   Bessel sequence,
 tight   frame) for $\mathcal{U}$.

Next we assume that
$\{\Lambda_j:\,j\in\bbJ\}$ is a g-Riesz basis for $\mathcal{U}$.
Since $\{e_{j,k}:\,k\in\bbK_j\}$ is an orthonormal basis for $\mathcal{V}_j$,
every $g_j\in\mathcal{V}_j$ has an expansion of the form
$g_j = \sum_{k\in\bbK_j} c_{j,k}e_{j,k}$,
where $\{c_{j,k}:\, k\in\bbK_j\}\in\ell^2(\bbK_j)$. It follows that
\[
  A\sum_{j\in \bbJ_1}\|g_j\|^2
   \le \left\| \sum_{j\in\bbJ_1}  \Lambda_j^* g_j\right\|^2
 \le   B\sum_{j\in \bbJ_1}\|g_j\|^2
\]
is equivalent to
\[
  A\sum_{j\in \bbJ_1}\sum_{k\in\bbK_j} |c_{j,k}|^2
   \le \left\| \sum_{j\in\bbJ_1} \sum_{k\in\bbK_j}  c_{j,k}u_{j,k}\right\|^2
 \le   B\sum_{j\in \bbJ_1}\sum_{k\in\bbK_j} |c_{j,k}|^2.
\]
On the other hand,
we see from $\Lambda_j f =\sum_{k\in\bbK_j} \langle f, u_{j,k}\rangle e_{j,k}$
that
$\{f:\, \Lambda_jf=0, j\in\bbJ\}
=\{f:\, \langle f, u_{j,k}\rangle=0, j\in\bbJ, k\in\bbK_j\}$.
Hence
$\{\Lambda_j:\,j\in\bbJ\}$ is g-complete if and only if $\{u_{j,k}:\,j\in\bbJ,k\in\bbK_j\}$
is   complete.
Therefore,
$\{\Lambda_j:\,j\in\bbJ\}$ is a g-Riesz basis  if and only if $\{u_{j,k}:\,j\in\bbJ,k\in\bbK_j\}$
is a    Riesz basis.

Now  we assume that
$\{\Lambda_j:\,j\in\bbJ\}$ is a g-orthonormal basis.
It follows from (\ref{eq:orth:a}) and (\ref{eq:ujk}) that
\begin{eqnarray*}
  \langle u_{j_1,k_1}, u_{j_2,k_2}\rangle
  &=&\langle \Lambda_{j_2} u_{j_1,k_1}, e_{j_2,k_2}\rangle \\
  &=&\overline{\langle   \Lambda_{j_2}^*e_{j_2,k_2}, u_{j_1,k_1}\rangle} \\
  &=&\overline{\langle \Lambda_{j_1} \Lambda_{j_2}^* e_{j_2,k_2}, e_{j_1,k_1}\rangle} \\
  &=&\langle \Lambda_{j_1}^* e_{j_1,k_1}, \Lambda_{j_2}^* e_{j_2,k_2}\rangle \\
   &=& \delta_{j_1,j_2}\delta_{k_1,k_2}, \qquad \forall j_1,j_2\in\bbJ,\,\,
    k_1\in\bbK_{j_1}, k_2\in\bbK_{j_2}.
\end{eqnarray*}
Hence
$\{u_{j,k}:\,j\in\bbJ, k\in\bbK_j\}$ is an  orthonormal sequence.
Moreover, observe  that
\[
 \|f\|^2= \sum_{j\in\bbJ} \|\Lambda_j f\|^2 = \sum_{j\in\bbJ} \sum_{k\in\bbK_j}
          |\langle  f, u_{j,k}\rangle|^2, \qquad \forall f\in\mathcal{U}.
\]
We have
$\{u_{j,k}:\,j\in\bbJ, k\in\bbK_j\}$ is an orthonormal basis.

For the converse, we need only to show that (\ref{eq:orth:a}) holds.
In fact, we see from (\ref{eq:rep:phi*}) that
for any $j_1\ne j_2\in \bbJ$,
$g_{j_1}\in\mathcal{V}_{j_1}$
and $g_{j_2}\in\mathcal{V}_{j_2}$,
\[
\langle \Lambda_{j_1}^* g_{j_1}, \Lambda_{j_2}^* g_{j_2}\rangle
= \left\langle \sum_{k_1\in\bbK_{j_1}} \langle g_{j_1}, e_{j_1,k_1}\rangle u_{j_1,k_1},
   \sum_{k_2\in\bbK_{j_2}} \langle g_{j_2}, e_{j_2,k_2}\rangle u_{j_2,k_2} \right\rangle
   =0
\]
and for $g_1, g_2\in \mathcal{V}_j$,
\[
\langle \Lambda_{j}^* g_{1}, \Lambda_{j}^* g_{2}\rangle
= \left\langle \sum_{k_1\in\bbK_{j}} \langle g_{1}, e_{j,k_1}\rangle u_{j,k_1},
   \sum_{k_2\in\bbK_{j}} \langle g_{2}, e_{j,k_2}\rangle u_{j,k_2} \right\rangle
= \langle g_1, g_2\rangle.
\]
Now the conclusion follows.

(ii).\,\,Since the cardinity
of a frame is no less than that of a basis, we
have
$\#\{u_{j,k}:\,j\in\bbJ,k\in\bbK_j\} \ge \dim \mathcal{U}$.
Hence
$       \sum_{j\in\bbJ} \dim \mathcal{V}_j \ge \dim \mathcal{U}$.
Moreover, we see from (i) that the equality holds whenever
$\{\Lambda_j:\,j\in\bbJ\}$ is a g-Riesz basis.

(iii).\,\, We see from (\ref{eq:rep:phij})
and (\ref{eq:rep:phi*}) that
\begin{eqnarray*}
 \sum_{j\in\bbJ} \Lambda_j^* \Lambda_j f
 &=& \sum_{j\in\bbJ} \sum_{k\in\bbK_j} \langle \Lambda_j f,  e_{j,k} \rangle u_{j,k}\\
 &=& \sum_{j\in\bbJ} \sum_{k\in\bbK_j} \left\langle
        \sum_{k'\in\bbK_j} \langle f, u_{j,k'}\rangle e_{j,k'},  \,\ e_{j,k}
     \right\rangle u_{j,k}\\
 &=& \sum_{j\in\bbJ} \sum_{k\in\bbK_j}
         \langle f, u_{j,k}\rangle u_{j,k}, \qquad \forall f\in\mathcal{U}.
\end{eqnarray*}
Hence the g-frame operator for $\{\Lambda_j:\,j\in\bbJ\}$
coincides with the   frame operator for
$\{u_{j,k}:\,j\in\bbJ,k\in\bbK\}$.

(iv).\,\, This  is a consequence of (i) and (iii).
The proof is over.
\end{proof}

The followings are  immediate consequences. We leave the proofs to interested readers.
\begin{Corollary}
\label{co:bessel}
$\{\Lambda_j:\,j\in\bbJ\}$ is
a g-Bessel sequence with an upper bound $B$ if and only if
for any finite subset $\bbJ_1\subset \bbJ$,
\[
      \bigg\|\sum_{j\in\bbJ_1} \Lambda_j^* g_j\bigg\|^2
      \le B \sum_{j\in\bbJ_1} \|g_j\|^2,
       \qquad
     g_j\in \mathcal{V}_j.
\]
\end{Corollary}

\begin{Corollary}
A g-Riesz basis $\{\Lambda_j:\,j\in\bbJ\}$ is
an exact g-frame.
Moreover,
it is
 g-biorthonormal with respect to its
dual $\{\tilde\Lambda_j:\,j\in\bbJ\}$ in the following sense
\[
 \langle \Lambda_{j_1}^* g_{j_1}, \tilde\Lambda_{j_2}^* g_{j_2}\rangle =
    \delta_{j_1, j_2}  \langle g_{j_1}, g_{j_2}\rangle ,
      \qquad
      \forall  j_1, j_2\in\bbJ, \quad g_{j_1}\in\mathcal{V}_{j_1},
            \quad g_{j_2}\in\mathcal{V}_{j_2}.
\]
\end{Corollary}

\begin{Corollary}\label{Co:g-riesz:a}
A sequence $\{\Lambda_j:\,j\in\bbJ\}$
is a g-Riesz basis for $\mathcal{U}$
with respect to $\seqV$
if and only if there is
a g-orthonormal basis $\{Q_j:\,j\in\bbJ\}$ for $\mathcal{U}$
and
a bounded invertible linear operator $T$ on $\mathcal{U}$
such that
\(
       \Lambda_j = Q_j T, j\in\bbJ.
\)
\end{Corollary}

\begin{proof}\,\,
Let $\{e_{j,k}:\,k\in\bbK_j\}$ be an orthonormal basis for $\mathcal{V}_j$, $j\in\bbJ$.
First, we assume that
$\{\Lambda_j:\,j\in\bbJ\}$
is a g-Riesz basis for $\mathcal{U}$.
By Theorem \ref{thm:equivalence}, we can find some   Riesz basis $\{u_{j,k}:\,j\in\bbJ, k\in\bbK_j\}$
for $\mathcal{U}$
such that
\[
     \Lambda_j f = \sum_{k\in\bbK_j} \langle f, u_{j,k}\rangle e_{j,k}.
\]
Take an    orthonormal basis
$\{u^{\circ}_{j,k}:\,j\in\bbJ, k\in\bbK_j\}$ for $\mathcal{U}$
 and define the operator $T$ on $\mathcal{U}$
by
\[
    T^*  u^{\circ}_{j,k} = u_{j,k}.
\]
Obviously, $T$ is a bounded invertible operator.
Let $Q_j\in \mathcal{L}(\mathcal{U}, \mathcal{V}_j)$ be such that
$Q_j f = \sum_{k\in\bbK_j} \langle f, u^{\circ}_{j,k}\rangle e_{j,k}$.
By Theorem \ref{thm:equivalence},
$\{Q_j:\,j\in\bbJ\}$ is a g-orthonormal basis.
Moreover,
for any $f\in\mathcal{U}$,
\[
       Q_j Tf =
           \sum_{k\in\bbK_j} \langle Tf, u^{\circ}_{j,k}\rangle e_{j,k}
       = \sum_{k\in\bbK_j} \langle f, T^*u^{\circ}_{j,k}\rangle e_{j,k}
       = \sum_{k\in\bbK_j} \langle f, u_{j,k}\rangle e_{j,k}
       =\Lambda_j f.
\]
Hence $\Lambda_j = Q_j T, \forall j\in\bbJ$.

Next we assume that
$\{Q_j:\,j\in\bbJ\}$ is a g-orthonormal basis
and $\Lambda_j = Q_j T$ for some
  bounded invertible operator $T$.
Then
$\{\Lambda_j:\,j\in\bbJ\}$ is g-complete in $\mathcal{U}$
and we can find some   orthonormal basis $\{u^{\circ}_{j,k}:\,j\in\bbJ,
k\in\bbK_j\}$
for $\mathcal{U}$ such that
$Q_j f =  \sum_{k\in\bbK_j} \langle f,  u^{\circ}_{j,k}\rangle e_{j,k}$.
Hence
$\Lambda_j f = \sum_{k\in\bbK_j} \langle Tf,  u^{\circ}_{j,k}\rangle e_{j,k}
=\sum_{k\in\bbK_j} \langle f, T^* u^{\circ}_{j,k}\rangle e_{j,k}$.
Now we see from Theorem \ref{thm:equivalence} that
$\{\Lambda_j:\,j\in\bbJ\}$ is a g-Riesz basis.
\end{proof}

\subsection{Excess of g-frames}
By Theorem \ref{thm:equivalence},
g-frames, g-Riesz bases and g-orthonormal bases
have similar properties
as frames, Riesz bases and orthonormal bases, respectively.
However, not all the properties are similar.
For example,   Riesz bases are equivalent to exact   frames.
But it is not the case for g-Riesz bases and exact g-frames.
In fact, we see from Theorems \ref{thm:equivalence} that  a g-Riesz basis is also an exact g-frame
while the converse is not true, which is not surprising since one element of a g-frame
might correspond to  several elements of the induced   frame.

\begin{Example}
Let $\{\varphi_j:\,j\in\bbJ\}$ be a Riesz basis for some Hilbert space $\mathcal{H}$.
Define $\Lambda_j:\, \mathcal{H}\mapsto \mathbb C^2$ as follows:
\[
  \Lambda_j f = (\langle f, \varphi_j\rangle, 0)^T.
\]
Then $\{\Lambda_j:\,j\in\bbJ\}$
is an exact g-frame. By Theorem \ref{thm:equivalence}, it is not a g-Riesz basis for
$\mathcal{H}$ with respect to $\mathbb C^2$.
However, it is a g-Riesz basis for
$\mathcal{H}$ with respect to
$\mathbb{C}\times \{0\}$.
\end{Example}

The above example shows that an exact g-frame may be a g-Riesz basis when we change the reference.
Does this hold in general? The answer is negative.

\begin{Example}
Let $\{\varphi_j:\,j\in\bbZ\}$ a Riesz basis for some Hilbert space $\mathcal{H}$.
Define $\Lambda_j:\, \mathcal{H}\mapsto \mathbb{C}^3$ as follows:
\[
  \Lambda_j f = (\langle f, \varphi_{2j-1}\rangle, \langle f, \varphi_{2j}\rangle,
      \langle f, \varphi_{2j+1}\rangle)^T.
\]
Then $\{\Lambda_j:\,j\in\bbZ\}$
is an exact g-frame.
However,
$\{\Lambda_j:\,j\in\bbZ\}$ is not a g-Riesz basis for
$\mathcal{H}$ with respect to any $\seqV$,
thanks to Theorem \ref{thm:equivalence}.
\end{Example}

On the other hand, it is well known (e.g., see \cite{Y}) that
a    frame   either remains a    frame or is incomplete whenever any one of its elements is removed.
It is neither the case for g-frames due to the same reason.
The following is a counterexample.
\begin{Example}
\label{ex:ex21}
Let $g(x) = e^{-x^2/2}$ be the Gaussian and
$\{\alpha_{m,n}:\,m,n\in\bbZ\}$ be an orthonormal basis for $\ell^2(\bbZ^2)$.
Define
\begin{eqnarray*}
 \Lambda_j f&=& \sum_{m,n\in\bbZ} \langle f(x), e^{i2\pi m x} g(x-2n-j)\rangle \alpha_{m,n},
     \quad j=1,2, \\
 \Lambda_3 f &=& \sum_{m,n\in\bbZ} \langle f(x), e^{i2\pi m x} g(x-n+1/2)\rangle \alpha_{m,n},
  \qquad
      f\in \LTR.
\end{eqnarray*}
We see from
Theorem \ref{thm:equivalence} and the frame theory (e.g., see \cite[p. 84--86]{D2}) that
  $\{\Lambda_1, \Lambda_2, \Lambda_3\}$
is a    g-frame for $\LTR$ with respect to $\ell^2(\bbZ^2)$.
However,  $\{\Lambda_1, \Lambda_2\}$  is not a g-frame but g-complete.
\end{Example}

A natural problem arises: in which case a subsequence
of a g-frame for which only one element is removed
is a g-frame  or not?
To this problem, we have the following.

\begin{Theorem}
\label{thm:frame2}
Let
$\{\Lambda_j:\,j\in\bbJ\}$ be a g-frame for $\mathcal{U}$ with respect to
$\seqV$ and
$\{\tilde\Lambda_j:\,j\in\bbJ\}$ be  the canonical dual g-frame.
Suppose that $j_0\in\bbJ$.

\begin{enumerate}
\item If there is some $g_0\in \mathcal{V}_{j_0}\setminus\{0\}$ such that
$ \tilde{\Lambda}_{j_0}  \Lambda_{j_0}^* g_0=g_0$,
then $\{\Lambda_j:\,j\in\bbJ,j\ne j_0\}$ is   not g-complete in $\mathcal{U}$.

\item If there is some $f_0\in \mathcal{U}\setminus\{0\}$ such that
$ {\Lambda}_{j_0}^*  \tilde\Lambda_{j_0} f_0 = f_0$,
then $\{\Lambda_j:\,j\in\bbJ,j\ne j_0\}$ is   not g-complete in $\mathcal{U}$.

\item If $I - \Lambda_{j_0}\tilde{\Lambda}_{j_0}^*$ or $
I - \tilde\Lambda_{j_0}{\Lambda}_{j_0}^*$
 is bounded invertible on $\mathcal{V}_{j_0}$,
then $\{\Lambda_j:\,j\in\bbJ,j\ne j_0\}$ is   a g-frame for $\mathcal{U}$.
\end{enumerate}
\end{Theorem}

\begin{proof}\,\,
(i).\,\,
Since $\Lambda_{j_0}^* g_0\in\mathcal{U}$, we have
\[
     \Lambda_{j_0}^* g_0 = \sum_{j\in\bbJ} \Lambda_j^* \tilde{\Lambda}_j
      \Lambda_{j_0}^* g_0
\]
Hence,
\(
   0 = \sum_{j\in\bbJ, j\ne j_0}\Lambda_j^* \tilde{\Lambda}_j  \Lambda_{j_0}^* g_0.
\)
Put $v_{j_0,j} = \delta_{j_0,j} g_0$. We have
\[
     \Lambda_{j_0}^* g_0 = \sum_{j\in\bbJ} \Lambda_j^* v_{j_0,j}.
\]
It follows from Lemma \ref{Lm:g-frame1} that
\[
  \sum_{j\in\bbJ} \|v_{j_0,j}\|^2
   = \sum_{j\in\bbJ} \|\tilde{\Lambda}_j  \Lambda_{j_0}^* g_0\|^2
   + \sum_{j\in\bbJ} \|\tilde{\Lambda}_j  \Lambda_{j_0}^* g_0 - v_{j_0,j}\|^2
\]
Consequently,
\[
 \|g_0\|^2
   = \|g_0\|^2
   + 2 \sum_{j\ne j_0} \|\tilde{\Lambda}_j  \Lambda_{j_0}^* g_0 \|^2
\]
Hence, $\tilde{\Lambda}_j  \Lambda_{j_0}^* g_0 = 0$.
Therefore,
\(
  {\Lambda}_j  \tilde \Lambda_{j_0}^* g_0
  = {\Lambda}_j  S^{-1} \Lambda_{j_0}^* g_0
  = \tilde{\Lambda}_j  \Lambda_{j_0}^* g_0
  =0\),
  $j\ne j_0$.
But
$\langle \Lambda_{j_0}^* g_0, \tilde \Lambda_{j_0}^* g_0\rangle
=\langle \tilde{\Lambda}_{j_0} \Lambda_{j_0}^* g_0,  g_0\rangle
=\|g_0\|^2>0$, which implies that
$\tilde \Lambda_{j_0}^* g_0\ne 0$.
Hence
$\{\Lambda_j:\,j\in\bbJ,j\ne j_0\}$ is not g-complete in $\mathcal{U}$.

(ii) Since
$ {\Lambda}_{j_0}^*  \tilde\Lambda_{j_0} f_0 = f_0 \ne 0$,
we have $\tilde\Lambda_{j_0} f_0\ne 0$
and
$ \tilde\Lambda_{j_0}{\Lambda}_{j_0}^*  \tilde\Lambda_{j_0} f_0 = \tilde\Lambda_{j_0}f_0$.
Now the conclusion follows from (i).

(iii).
Since $\tilde\Lambda_j = \Lambda_j S^{-1}$, where $S$ is the g-frame operator for
$\{\Lambda_j:\,j\in\bbJ\}$,
we have
\[
     I - \Lambda_{j_0} \tilde{\Lambda}_{j_0}^*
=  I - \Lambda_{j_0} S^{-1} {\Lambda}_{j_0}^*
=  I - \tilde\Lambda_{j_0} {\Lambda}_{j_0}^*.
\]
Let $A$ and $B$ be the lower and upper  frame bounds for
$\{\Lambda_j:\,j\in\bbJ\}$, respectively.
For any $f\in\mathcal{U}$, we have
\[
       f = \sum_{j\in\bbJ} \tilde{\Lambda}_{j}^* \Lambda_j f.
\]
Hence
\[
      \Lambda_{j_0} f = \sum_{j\in\bbJ}  \Lambda_{j_0}\tilde{\Lambda}_{j}^* \Lambda_j f.
\]
Therefore,
\begin{equation}\label{eq:t3}
      (I - \Lambda_{j_0} \tilde{\Lambda}_{j_0}^*  )\Lambda_{j_0} f
       = \sum_{j\ne j_0}  \Lambda_{j_0}\tilde{\Lambda}_{j}^* \Lambda_j f.
\end{equation}
Note that
\begin{eqnarray*}
\left\|\sum_{j\ne j_0}  \Lambda_{j_0}\tilde{\Lambda}_{j}^* \Lambda_j f\right\|^2
&=&\sup_{g\in \mathcal{V}_{j_0}, \|g\|=1}
    \left|\left\langle\sum_{j\ne j_0}  \Lambda_{j_0}\tilde{\Lambda}_{j}^* \Lambda_j f,
    g\right\rangle\right|^2\\
&=&\sup_{\|g\|=1}
    \left|\sum_{j\ne j_0} \left\langle  \Lambda_j f,
    \tilde{\Lambda}_{j} \Lambda_{j_0}^*g\right\rangle\right|^2\\
&\le&
\sum_{j\ne j_0}  \| \Lambda_j f\|^2\cdot
\sup_{\|g\|=1}
 \sum_{j\in\bbJ}   \|    \tilde{\Lambda}_{j} \Lambda_{j_0}^*g\|^2\\
&\le&
   \frac{1}{A} \|\Lambda_{j_0}\|^2 \sum_{j\ne j_0}  \| \Lambda_j f\|^2.
\end{eqnarray*}
We see from (\ref{eq:t3}) that
\[
    \| \Lambda_{j_0} f\|^2
    \le      \|(I - \Lambda_{j_0}\tilde{\Lambda}_{j_0}^*  )^{-1}\|^2
       \frac{1}{A} \|\Lambda_{j_0}\|^2 \sum_{j\ne j_0}  \| \Lambda_j f\|^2.
\]
Hence
\[
 \sum_{j\in \bbJ}  \| \Lambda_j f\|^2 \le C \sum_{j\ne j_0}  \| \Lambda_j f\|^2.
\]
Therefore,
\[
 \frac{A}{C} \|f\|^2 \le \sum_{j\ne j_0}  \| \Lambda_j f\|^2 \le  B\|f\|^2, \qquad \forall f\in\mathcal{U}.
\]
This completes the proof.
\end{proof}

\begin{Corollary}
\label{Co:frame2}
Let
$\{\Lambda_j:\,j\in\bbJ\}$ be a g-frame for $\mathcal{U}$ with respect to
$\seqV$.
If $\dim \mathcal{V}_j < +\infty$, $j\in\bbJ$,
then
$\{\Lambda_j:\,j\in\bbJ,j\ne j_0\}$ is  either g-incomplete in $\mathcal{U}$
or a g-frame for $\mathcal{U}$
 for any $j_0\in\bbJ$.
\end{Corollary}

\begin{proof}
If there is some $g_0\in \mathcal{V}_{j_0}\setminus\{0\}$ such that
$ \tilde{\Lambda}_{j_0}  \Lambda_{j_0}^* g_0=g_0$,
then Theorem \ref{thm:frame2} (i) shows that
$\{\Lambda_j:\,j\in\bbJ,j\ne j_0\}$ is   not g-complete in $\mathcal{U}$.
Otherwise,  $I - \tilde{\Lambda}_{j_0}  \Lambda_{j_0}^*$ is injective. Consequently,
$(I - \Lambda_{j_0} \tilde{\Lambda}_{j_0}^*)\mathcal{V}_{j_0}$ is dense
in $\mathcal{V}_{j_0}$.
Since $\dim \mathcal{V}_{j_0}<+\infty$,
we have
$(I - \Lambda_{j_0} \tilde{\Lambda}_{j_0}^*)\mathcal{V}_{j_0}=\mathcal{V}_{j_0}$.
Therefore,
$I - \Lambda_{j_0} \tilde{\Lambda}_{j_0}^*$ is bounded invertible.
Now the conclusion follows from
Theorem \ref{thm:frame2} (iii).
\end{proof}

\subsection{Equivalence between
stable space splittings and g-frames}

Stable space splittings are  generalizations of frames
which lead to a better understanding of iterative solvers
(multigrid/multilevel resp. domain decomposition methods)
for large-scale discretization of elliptic operator equations
(see \cite{Osw2} and references therein).
Here we prove that
stable space splittings  are equivalent to g-frames.

Let $\mathcal{V}$ and $\mathcal{V}_j, j\in\bbJ$ be  Hilbert spaces.
Let $b_j$ be a bilinear form on $\mathcal{V}_j\times \mathcal{V}_j$
satisfying
\begin{equation}\label{eq:bj}
     b_j(u,u) \ge C_j \|u\|^2\, \textrm{and}\, b_j(u,v)=b_j(v,u) \le C_j' \|u\|\cdot\|v\|,
     \,\forall u,v\in \mathcal{V}_j.
\end{equation}
Suppose that $R_j\in\mathcal{L}(\mathcal{V}_j, \mathcal{V})$.
Recall that a
system $\{(\{\mathcal{V}_j, b_j\},R_j):\,j\in\bbJ\}$ is called a stable space splitting of
$\mathcal{V}$ if there are some positive constants $C, C'$ such that
\[
     C \|u\|^2 \le \inf_{u=\sum_{j\in\bbJ} R_j u_j} \sum_{j\in\bbJ} b_j(u_j, u_j)
      \le C'\|u\|^2, \qquad \forall u\in \mathcal{V},
\]
where $u_j\in \mathcal{V}_j$.
It was shown in \cite[Theorem 4]{Osw2} (see also \cite[Pages 73-75]{Osw1})
that a stable space splitting $\{(\{\mathcal{V}_j, b_j\},R_j):\,j\in\bbJ\}$ satisfies
$A \le \sum_{j\in\bbJ} R_jR_j^* \le B$ for some constants $A,B>0$, which is equivalent to
\[
   A\|u\|^2 \le \sum_{j\in\bbJ} \| R_j^* u\|^2
   \le B\|u\|^2, \qquad \forall u\in \mathcal{V}.
\]
Hence $\{R_j^*:\,j\in\bbJ\}$ is a g-frame for
 $\mathcal{V}$ with respect to $\seqV$.

For the converse, we need
$b_j(u,u)$ to be uniformly bounded, i.e., we assume that
\begin{equation}\label{eq:bjs}
     C_1 \|u\|^2 \le b_j(u,u) \le C_2 \|u\|^2, \qquad \forall u\in\mathcal{V}.
\end{equation}
Suppose that $\{\Lambda_j:\,j\in\bbJ\}$ is a
g-frame for
 $\mathcal{V}$ with respect to $\seqV$.
Let $A$ and $B$ be the frame bounds  and $\{\tilde\Lambda_j:\,j\in\bbJ\}$
be the canonical dual g-frame.
Put $R_j = \Lambda_j^*$.
We see from (\ref{eq:bjs})  and Lemma \ref{Lm:g-frame1} that
\begin{eqnarray*}
 \inf_{u=\sum_{j\in\bbJ} R_j u_j} \sum_{j\in\bbJ} b_j(u_j, u_j)
 &\ge&  \inf_{u=\sum_{j\in\bbJ} R_j u_j} \sum_{j\in\bbJ} C_1\|u_j\|^2 \\
 &=& \sum_{j\in\bbJ} C_1\|\tilde\Lambda_j u\|^2
\ge \frac{C_1}{B} \|u\|^2.
\end{eqnarray*}
Similarly we can prove that
\[
 \inf_{u=\sum_{j\in\bbJ} R_j u_j} \sum_{j\in\bbJ} b_j(u_j, u_j)
\le\frac{C_2}{A} \|u\|^2.
\]
Hence
$\{(\{\mathcal{V}_j, b_j\},R_j):\,j\in\bbJ\}$ is a stable space splitting.

\section{Applications of g-frames}
\subsection{Atomic resolution  of bounded linear operators}

Here we give an application of g-frames.

Let $\{\Lambda_j:\,j\in\bbJ\}$ be a g-frame for $\mathcal{U}$ with respect to $\seqV$.
Suppose that
$\{\tilde{\Lambda}_j:\,j\in\bbJ\}$ is the canonical dual g-frame.
Then for any $f\in \mathcal{U}$, we have
\[
   f = \sum_{j\in\bbJ} \Lambda_j^* \tilde\Lambda_j f
     = \sum_{j\in\bbJ} \tilde\Lambda_j^* \Lambda_j f, \qquad f\in\mathcal{U}.
\]
It follows that
\begin{equation}\label{eq:Id}
    I_{\mathcal{U}} =
  \sum_{j\in\bbJ} \Lambda_j^* \tilde\Lambda_j
       = \sum_{j\in\bbJ} \tilde\Lambda_j^* \Lambda_j,
\end{equation}
where the convergence is in weak* sense.
Let $T$ be a bounded linear operator on $\mathcal{U}$.
We see from (\ref{eq:Id}) that
\begin{equation}\label{eq:res}
   T  = \sum_{j\in\bbJ} T\Lambda_j^* \tilde\Lambda_j
      = \sum_{j\in\bbJ} T\tilde\Lambda_j^* \Lambda_j
      = \sum_{j\in\bbJ} \Lambda_j^* \tilde\Lambda_j T
      = \sum_{j\in\bbJ} \tilde\Lambda_j^* \Lambda_j T.
\end{equation}
We call
(\ref{eq:res})  atomic resolutions of an operator $T$.

\subsection{Construction of   frames via g-frames}

Let $\{\Lambda_j:\,j\in\bbJ\}$ be a g-frame for $\mathcal{U}$
with respect to $\seqV$. We see from Theorem \ref{thm:equivalence}
that $\{u_{j,k}:\,j\in\bbJ,k\in\bbK_j\}
=\{\Lambda_j^* e_{j,k}:\,j\in\bbJ,k\in\bbK_j\}
$ is a frame for $\mathcal{U}$, where
$\{e_{j,k}:\,k\in\bbK_j\}$
is an orthonormal basis for $\mathcal{V}_j$.
However, it might be difficult to find
an orthonormal basis for $\mathcal{V}_j$ in practice.
Fortunately, the orthonormality is not necessary to get a frame.
In fact, we have the following.

\begin{Theorem}
\label{thm:const:frame}
Let $\{\Lambda_j:\,j\in\bbJ\}$
and $\{\tilde{\Lambda}_j:\,j\in\bbJ\}$
be a pair of dual g-frames for $\mathcal{U}$ with respect to $\seqV$
and
 $\{g_{j,k}:\,k\in\bbK_j\}$  and $\{\tilde g_{j,k}:\,k\in\bbK_j\}$
be a  pair of dual    frames for $\mathcal{V}_j$, respectively.
Then
$\{\Lambda_j^*g_{j,k}:\,j\in\bbJ,k\in\bbK_j\}$
and
$\{\tilde\Lambda_j^*\tilde g_{j,k}:\,j\in\bbJ,k\in\bbK_j\}$
are a pair of dual    frames for $\mathcal{U}$.

Moreover, suppose that
$\{\Lambda_j:\,j\in\bbJ\}$
and $\{\tilde{\Lambda}_j:\,j\in\bbJ\}$ are canonical dual g-frames,
$\{g_{j,k}:\,k\in\bbK_j\}$  and $\{\tilde g_{j,k}:\,k\in\bbK_j\}$
are canonical dual frames,
and that
$\{g_{j,k}:\,k\in\bbK_j\}$
is a tight g-frame with frame bounds $A_j=B_j=A, \forall j\in\bbJ$.
Then
$\{\Lambda_j^*g_{j,k}:\,j\in\bbJ,k\in\bbK_j\}$
and
$\{\tilde\Lambda_j^*\tilde g_{j,k}:\,j\in\bbJ,k\in\bbK_j\}$
are canonical dual  frames.
\end{Theorem}

\begin{proof} Note that
\[
     \langle f, \Lambda_j^*g_{j,k}\rangle     =      \langle \Lambda_j f, g_{j,k}\rangle.
\]
It is easy to see  that both
$\{\Lambda_j^*g_{j,k}:\,j\in\bbJ,k\in\bbK_j\}$
and
$\{\tilde\Lambda_j^*\tilde g_{j,k}:\,j\in\bbJ,k\in\bbK_j\}$
are   frames for $\mathcal{U}$.
On the other hand,
 For any $f\in\mathcal{U}$, we have
\[
   \sum_{j\in\bbJ}\sum_{k\in\bbK_j}
      \langle f, \Lambda_j^*g_{j,k}\rangle \tilde \Lambda_j^* \tilde g_{j,k}
 =   \sum_{j\in\bbJ}
      \tilde \Lambda_j^* \sum_{k\in\bbK_j}
      \langle\Lambda_j f, g_{j,k}\rangle  \tilde g_{j,k} \\
 =  \sum_{j\in\bbJ}
      \tilde \Lambda_j^*
       \Lambda_j f\\
 =    f.
\]
Similarly we can get that
\[
   \sum_{j\in\bbJ}\sum_{k\in\bbK_j}
      \langle f, \tilde \Lambda_j^* \tilde g_{j,k}\rangle \Lambda_j^*g_{j,k}
    =f.
\]
Hence
$\{\Lambda_j^*g_{j,k}:\,j\in\bbJ,k\in\bbK_j\}$
and
$\{\tilde\Lambda_j^*\tilde g_{j,k}:\,j\in\bbJ,k\in\bbK_j\}$
are dual   frames for $\mathcal{U}$.

Next we assume that
$\{\Lambda_j:\,j\in\bbJ\}$
and $\{\tilde{\Lambda}_j:\,j\in\bbJ\}$
are  canonical dual g-frames
and $\{g_{j,k}:\,k\in\bbK_j\}$ is a tight   frame with frame bounds $A_j=B_j=A$, $j\in\bbJ$.
Then $\tilde{g}_{j,k} = \frac{1}{A} g_{j,k}$.
Let $S_{\Lambda}$ and $S_{\Lambda,g}$ be the frame operators associated with
$\{\Lambda_j:\,j\in\bbJ\}$
and
 $\{\Lambda_j^* g_{j,k}:\,j\in\bbJ, k\in\bbK_j\}$, respectively.
Then we have
\begin{eqnarray*}
   S_{\Lambda,g} f
&=&
   \sum_{j\in\bbJ}\sum_{k\in\bbK_j}
      \langle f, \Lambda_j^* g_{j,k}\rangle   \Lambda_j^*  g_{j,k} \\
&=&\sum_{j\in\bbJ}\Lambda_j^*  \sum_{k\in\bbK_j}
      \langle \Lambda_j f,  g_{j,k}\rangle    g_{j,k}\\
&=& A \sum_{j\in\bbJ}\Lambda_j^*   \Lambda_j f\\
&=& A S_{\Lambda}f, \qquad \forall f\in\mathcal{U}.
\end{eqnarray*}
Hence
\[
 S_{\Lambda,g}^{-1}\Lambda_j^* g_{j,k}
 = \frac{1}{A}S_{\Lambda}^{-1}\Lambda_j^* g_{j,k}
 = \tilde\Lambda_j^* \tilde g_{j,k}, \qquad j\in\bbJ, k\in\bbK_j.
\]
This completes the proof.
\end{proof}

\section*{Acknowledgments}

Part  of this work was done while the  author was visiting
NuHAG at the Faculty of Mathematics, University of Vienna  and
the Erwin Schr\"odinger International Institute for Mathematical Physics (ESI),
Vienna.
He thanks NuHAG and ESI for hospitality
and support.
At ESI the author met Peter Oswald who introduced him
to the concept of stable space splittings.
He thanks Ole Christensen for  carefully reading the paper
and for relevant comments and suggestions.
Particular  thanks
go to Hans G.~Feichtinger for many helpful suggestions
and for introducing several other generalizations of frames.


\begin{thebibliography}{EE}

\bibitem{ACM}
A.~Aldroubi, C.~Cabrelli, and U.~Molter,
{Wavelets on irregular grids with arbitrary dilation matrices and
frame atomics for $L^2(\mathbb{R}^d)$},
\textit{Appl. Comput. Harmon. Anal.},
\textbf{17}(2004), 119--140.



\bibitem{AK}
M.S.~Asgari and A.~Khosravi,
Frames and bases of subspaces in Hilbert spaces,
\textit{J. Math. Anal. Appl.},
\textbf{308}(2005), 541--553.





\bibitem{CK}
P.~Casazza and G.~Kutyniok,
{Frames of subspaces},
in:
\textit{Wavelets, Frames and Operator Theory},
Vol. 345, American Mathematical Society, 2004, 87--113.


\bibitem{C6}O.~Christensen,
\textit{An Introduction to Frames and Riesz Bases},
Birkh\"auser, Boston, 2003.

\bibitem{CE}
O.~Christensen and Y.C.~Eldar,
Oblique dual frames and shift-invariant spaces,
\textit{Appl. Comput. Harmon. Anal.},
\textbf{17}(2004), 48--68.

\bibitem{CGY}
I.~Daubechies, A.~Grossmann, and Y.~Meyer,
{Painless nonorthogonal expansions},
\textit{J.Math. Phys.},
\textbf{27}(1986), 1271--1283.



\bibitem{D2}
I.~Daubechies,
\textit{Ten Lectures on Wavelets},
SIAM Philadelphia, 1992.

\bibitem{DFG}
M.~D\"orfler, H.G.~Feichtinger and K.~Gr\"ochenig,
Time-frequency partitions for the Gelfand triple,
preprint.

\bibitem{DS}
R.J.~Duffin and A.C.~Schaeffer,
{A class of nonharmonic Fourier series},
\textit{Trans. Amer. Math. Soc.},
\textbf{72}(1952), 341--366.



\bibitem{Eld}
Y.~Eldar,
{Sampling with arbitrary sampling and reconstruction spaces and
oblique dual frame vectors},
\textit{J. Fourier Anal. Appl.},
\textbf{9}(2003), 77--96.



\bibitem{FeiS}
H.G.~Feichtinger and T.~Strohmer, Eds.,
\textit{Gabor Analysis and Algorithms: Theory and Applications},
Birkh\"auser Inc., Boston, MA, 1998.

\bibitem{For1}
M.~Fornasier,
Quasi-orthogonal decompositions of structured frames,
\textit{J.  Math. Anal. Appl.},
\textbf{289}(2004), 180--199.

\bibitem{For2}
M.~Fornasier,
Decompositions of Hilbert spaces: local construction of global frames,
\textit{Proc. Int.
Conf. on Constructive function theory}, Varna (2002),
 B.~Bojanov Ed., DARBA, Sofia, 2003, 275--281.

\bibitem{Gr3} K.~Gr\"ochenig,
\textit{Foundations of Time-Frequency analysis},
Birk\"auser, Boston, 2001.



\bibitem{HL}
D.~Han, D.R.~Larson,
{Frames, bases and group representations},
\textit{Mem.  Amer.  Math.  Soc.},
\textbf{147}(2000), no.697.



\bibitem{HW}
C.~Heil and D.~Walnut,
{Continuous and discrete wavelet transforms},
\textit{SIAM Review},
\textbf{31}(1989), 628--666.


\bibitem{LO}
S.~Li and H.~Ogawa,
Pseudoframes for subspaces with applications,
\textit{J. Fourier Anal. Appl.},
\textbf{10}(2004), 409--431.


\bibitem{Osw1}
P.~Oswald,
\textit{Multilevel Finite Element Approximation: Theory
and Application},
Teubner Skripten zur Numerik,
Teubner, Stuttgart, 1994.


\bibitem{Osw2}
P.~Oswald,
{Frames and space splittings in Hilbert spaces},
preprint,
http://cm.bell-labs.com/who/poswald/bonn1.ps.gz






\bibitem{Y}
R.~Young,
\textit{An Introduction to Non-Harmonic Fourier Series,}
 Academic Press, New York, 1980.



\end{thebibliography}
\end{document}